\documentclass[12pt]{amsart}
\usepackage{amsmath, amscd}
\usepackage{amsfonts}
\begin{document}

\newcommand{\R}{\mathbb R}
\newcommand{\al}{\alpha}
\newcommand{\C}{\mathbb C}
\newcommand{\Po}{\mathbb P}
\newcommand{\N}{\mathbb N}
\newcommand{\ds}{\displaystyle}
\def\Re{{\sf Re}\,}
\def\Im{{\sf Im}\,}

\newtheorem{theorem}{Theorem}[section]

\theoremstyle{remark}
\newtheorem{remark}[theorem]{Remark}

\numberwithin{equation}{section}


\title[Quasi-parabolic one-resonant biholomorphisms]{Dynamics of quasi-parabolic one-resonant biholomorphisms}

\author{Filippo Bracci \and Feng Rong}

\address{Dipartimento di Matematica, Universit\`{a} di Roma``Tor Vergata", Via della Ricerca Scientifica 1, 00133 Roma, Italia}
\email{fbracci@mat.uniroma2.it}

\address{Department of Mathematics, Shanghai Jiao Tong University, 800 Dong Chuan Road, Shanghai, 200240, P.R. China}
\email{frong@sjtu.edu.cn}

\thanks{The first named author is partially supported by the ERC grant ``HEVO - Holomorphic Evolution Equations'' n. 277691. The second named author is partially supported by the National Natural Science Foundation of China (Grant No. 11001172), the Specialized Research Fund for the Doctoral Program of Higher Education of China (Grant No. 20100073120067) and the Scientific Research Starting Foundation for Returned Overseas Chinese Scholars.}

\begin{abstract}
In this paper we study the dynamics of germs of quasi-parabolic one-resonant  biholomorphisms of $\C^{n+1}$ fixing the origin, namely, those germs whose differential at the origin has one eigenvalue $1$ and the others having a one dimensional family of resonant relations. We define some invariants and give conditions which ensure the existence of attracting domains for such maps.
\end{abstract}

\maketitle

\section{Introduction}

Let $F$ be a germ of biholomorphism of $\C^{n+1}$ with a fixed point at the origin $O$. Assume that the eigenvalues of $dF_O$ are $\{1,\lambda_1,\cdots,\lambda_n\}$. The dynamics of such maps have been studied deeply by several authors.

In the tangent to the identity case, i.e. $\lambda_j=1$, $j=1,\ldots, n$, it has been proved by \'Ecalle \cite{Ec} and Hakim \cite{H:Parabolic} that generically there exist  ``parabolic curves'', namely, one-dimensional $F$-invariant analytic discs having the origin in their boundary and on which the dynamics is of parabolic type. Later, Abate \cite{Ab} (see also \cite{ABT}) proved that such parabolic curves always exist in dimension two, and Hakim \cite{H:Parabolic} gave also conditions which ensure the existence of basins of attraction of parabolic type (see also \cite{V} and \cite{L} for weaker conditions).

The semi-attractive case, namely when $|\lambda_j|<1$, $j=1,\ldots, n$, was studied by Fatou \cite{Fa}, Ueda \cite{U:Local1}, \cite{U:Local2}, Hakim \cite{H:Semi}, Rivi \cite{Ri:Semi} and the second named author \cite{R:Semi} (see also \cite{BZ:One}) who proved the existence of a basin of attraction. For further information we refer the reader to the survey papers \cite{B:Local} and \cite{A:Discrete}.

In case $|\lambda_j|=1$ and $\lambda_j$ not roots of unity for $j=1,\cdots,n$ the germ $F$ is a so-called \textit{quasi-parabolic} map. Quasi-parabolic maps  have been studied by the first named author and Molino  \cite{BM:Quasi} in dimension two and by the second named author in higher dimensions (\cite{R:Quasi1, R:Quasi2, R:Quasi3}). The focus of those papers was mainly in finding invariants which assure existence of parabolic curves.

In a different direction, let $G$ be a germ of holomorphic diffeomorphism of $\C^{n+1}$ fixing the origin whose differential $dG_0$ has eigenvalues $\mu_1,\cdots, \mu_{n+1}$ having the property that, for a fixed $m\leq n+1$, there exists a fixed multi-index $\al\in \N^{m}\times \{0\}^{n+1-m}$ for which the resonances $\mu_s=\mu^\beta:=\prod_{j=1}^{n+1} \mu_j^{\beta_j}$ for $1\leq s\leq m$ are precisely given by $\beta=k\al$ for some $k\geq 1$ arbitrary. Such a germ is called \textit{one-resonant with respect to $\{\mu_1,\cdots, \mu_m\}$}. In \cite{BZ:One}, the first named author with Zaitsev gave sharp invariants which assure the existence of basins of attraction of parabolic type for $G$ at $O$. Such a result has been extended to multi-resonant germs by the first named author with Raissy and Zaitsev in \cite{BRZ:Multi}.

In this paper we consider a mixed situation of quasi-parabolic and one-resonant, that is, we consider quasi-parabolic germs $F$ which are, in a certain sense, one-resonant with respect to $\{\lambda_1,\cdots, \lambda_n\}$. We assume that  there exists a fixed multi-index $\al\in \N^n$ such that if $\lambda_s=\lambda^\beta:=\prod_{j=1}^n \lambda_j^{\beta_j}$ for $1\leq s\leq n$, then necessarily $\beta=k\al$ for some $k\geq 1$ arbitrary. We call such a germ a \textit{quasi-parabolic one-resonant biholomorphism}.

The aim of this paper is to find invariants which ensure the existence of basins of attraction of parabolic type in such cases when the dynamics is somewhat mixed up. Since the eigenvalue $1$ enters in all resonances for the other eigenvalues, the setting is however different from the one-resonant case.

In order to state our result, we need to define some invariants. We refer to Section \ref{two} for details. Let $F$ be a germ of quasi-parabolic one-resonant biholomorphism fixing the origin $O$ and $\{1,\lambda_1,\cdots, \lambda_n\}$ the eigenvalues of $dF_O$. Let $2\leq \nu(F)<+\infty$ be the order of $F$ (if $\nu(F)=+\infty$ then $F$ has a holomorphic curve of fixed points passing through $O$, cf. \cite[Proposition 6.2]{BZ:One}). Suppose $F$ is dynamically separating with respect to the non-degenerate characteristic direction $[1:0\cdots:0]$. In \cite{BM:Quasi} and \cite{R:Quasi2} it has been proven that there exist at least $\nu(F)-1$ parabolic curves for $F$ tangent to $[1:0:\cdots:0]$ at $O$.

For $(z,w)\in \C \times \C^n$, we write $(z_1,w_{1,1},\cdots,w_{n,1}):=F(z,w)$. Under these conditions, one can perform  holomorphic changes of coordinates and finitely many blow-ups, in such a way that $F$ assumes the form
\begin{equation}\label{normal}\left\{
\begin{aligned}
z_1  = &z - \frac{1}{\nu(F)-1}z^{\nu(F)} + z^{\nu(F)+1}R_0(z,w),\\
w_{j,1} = &\lambda_j w_j - a_jz^{l_j}w^{\rho(F)\alpha}w_j - b_jz^{r_j}w_j \\&+\sum_{s>l_j, m>\rho(F)} p_{j,s}z^{s}w^{m\alpha}w_j+ \sum_{s>r_j}q_{j,s}z^{s}w_j+z^{\nu(F)+1}R_j(z,w),\ 1\le j\le n,
\end{aligned}
\right.
\end{equation}
where $(a_1,\cdots, a_n)\neq (0,\cdots, 0)$, $(b_1,\cdots, b_n)\neq (0,\cdots, 0)$, $\rho(F)\in \N$ is an invariant of $F$ which we call the {\sl one-resonant order of $F$}, $r_j\ge \nu(F)-1$ and the $R_j(z,w)$'s contain non-resonant terms of arbitrarily high order, for $j=1,\cdots, n$.

If $r_j=\nu(F)-1$ for some $j$, then the dynamics of the $(j+1)$-th component of the iterates of $F$, $w_{j,n}$, largely depends on $\Re(b_j)$ (cf. \cite{R:Quasi3}) and the one-resonant part does not play any role. Here instead we are interested in the degenerate situation, namely, when $r_j>\nu(F)-1$ for all $j=1,\cdots, n$.  If this is the case, we call $F$ a {\sl degenerately dynamically separating} quasi-parabolic germ. In such a case, the dynamics of $F$ relies also on the ``one-resonant'' part as we will show.

We assume also that $l_1=\cdots=l_n=:l(F)\leq \nu(F)-1$, and we call $l(F)$ the {\sl separation order}. In such a case, we let $A(F):=\sum_{j=1}^n \alpha_ja_j\lambda_j^{-1}$. Such a number is an invariant up to a non-zero scalar multiple, and $F$ is said to be \textit{non-degenerate with respect to $\{\lambda_1,\cdots,\lambda_n\}$} if $A(F)\neq 0$. Finally, if $A(F)\neq 0$, we say that $F$ is \textit{attracting with respect to $\{\lambda_1,\cdots,\lambda_n\}$} if $\Re(a_j\lambda_j^{-1}A(F)^{-1})>0$ for all $1\le j\le n$. Our main result is the following

\begin{theorem}\label{T:Main}
Let $F$ be a germ of quasi-parabolic one-resonant biholomorphism of $\C^{n+1}$ fixing $O$. Let $\{1,\lambda_1,\cdots,\lambda_n\}$ be the eigenvalues of  $dF_O$. Assume that $F$ is degenerately dynamically separating and attracting with respect to $\{\lambda_1,\cdots,\lambda_n\}$ and let $\rho(F)$ be its one-resonant order. Then $F$ has (at least) $\rho(F)$ disjoint basins of attraction at $O$.
\end{theorem}

\begin{remark}
More precisely, we get $(\nu(F)-1)\rho(F)$ disjoint basins of attraction when the separation order $l(F)=\nu(F)-1$ or 0, where $\nu(F)$ is the order of $F$, and $\rho(F)$ disjoint basins of attraction when $0<l(F)<\nu(F)-1$. See Remark \ref{R:l} for more details.
\end{remark}

\begin{remark}
One can easily extend the above result to biholomorphisms of $\C^{n+m+1}$, with eigenvalues $\{1,\lambda_1,\cdots,\lambda_n,\gamma_1,\cdots,\gamma_m\}$, where $\lambda_j$'s are as above and $|\gamma_i|<1$, $1\le i\le m$. This is similar to the easiest case for semi-attractive analytic transformations, which have been studied by several authors, as mentioned above.
\end{remark}

\medskip

Part of this work was done while the second named author was visiting IH\'{E}S and Dipartimento di Matematica, Universit\`{a} di Roma``Tor Vergata". He would like to thank the hosts for their hospitality and the institutes and K.C. Wong Education Foundation for the support.

\section{Quasi-parabolic one-resonant biholomorphisms}\label{two}

Let $F$ be a germ of quasi-parabolic one-resonant biholomorphism of $\C^{n+1}$ fixing the origin. Let $\{1,\lambda_1,\cdots, \lambda_n\}$ be the eigenvalues of $dF_0$. Note that  the one-resonant condition among $\{\lambda_1,\cdots,\lambda_n\}$ implies in particular that $\lambda_i\neq \lambda_j$ for $i\neq j$.

With a holomorphic change of coordinates we can assume that the non-resonant terms have order as high as we want. Thus, we can write $F$ as
$$\left\{
\begin{aligned}
z_1 & = z + f(z,w),\\
w_1 & = \Lambda w + g(z,w),
\end{aligned}
\right.$$
where $\Lambda=\textup{Diag}\{\lambda_1,\cdots,\lambda_n\}$, and $f(z,w)$ and $g(z,w)$ contain terms of order at least two.

Let $\nu$ ({\sl respectively} $\mu$) be the least of $i$ for terms $z^i$ in the expression of $z_1$ ({\sl resp.} $w_1$). If $\nu<\infty$ and $\mu\ge \nu$, then we say that $F$ is \textit{ultra-resonant}, and that the \textit{order} of $F$ is $\nu$. This is well-defined by \cite[Lemma 2.5]{BM:Quasi}, \cite[Lemma 2.3]{R:Quasi2}.

Suppose now that $F$ is ultra-resonant of order $\nu$. We assume that the vector $[v]=[1:0:\cdots:0]\in \Po^n$ is a \textit{non-degenerate characteristic direction} for $F$, {\sl i.e.} $F_\nu(v)=\lambda v$ for some $\lambda\neq 0$, where $F_\nu$ is the homogeneous part of $F$ of order $\nu$.

Write $w_1=(w_{1,1},\cdots,w_{n,1})$ and $w=(w_1,\cdots,w_n)$. We say that $F$ is \textit{dynamically separating} in the characteristic direction $[v]$ if there are no terms $z^iw_j$ with $i<\nu-1$ in the expression of $w_{j,1}$ for any $1\le j\le n$. This is well-defined by \cite[Lemma 2.10]{R:Quasi2}.

Now, after finitely many blow-ups centered at the point of the exceptional divisor given by the direction $[1:0:\cdots:0]$, we can bring $F$ to the form
\begin{equation}\label{k-inv}
\left\{
\begin{aligned}
z_1 & = z - \frac{1}{\nu-1}z^{\nu} + {\sf h.o.t.},\\
w_{j,1} & = \lambda_j w_j - a_jz^{l_j}w^{k\alpha}w_j - b_jz^{r_j}w_j + {\sf h.o.t.},\ \ \ 1\le j\le n,
\end{aligned}
\right.
\end{equation}
where ${\sf h.o.t.}$ means as usual {\sl higher order terms} and $r_j\ge \nu-1$ (since $F$ is dynamically separating).

The number $k\in \N$ is an invariant of $F$ under holomorphic changes of coordinates which preserve the form \eqref{k-inv} (cf. \cite[Remark 3.2]{BZ:One}). In fact, it is an invariant of $F$ even under blow-ups centered at the direction $[1:0:\cdots:0]$, because terms like $a_jz^{l_j}w^{k\alpha}w_j$ in $w_{j,1}$ are transformed into terms like $a_jz^{l_j+k|\al|}w^{k\alpha}w_j$ and the other terms either get higher degree or stay stable in case of terms like $z^{r_j}w_j$. In any case, no new terms like $z^{m}w^{k'\alpha}w_j$ may appear with $k'\leq k$.

We call the invariant $\rho(F):=k$ the  \textit{one resonant order of $F$}.

After performing a finite number of blow-ups, we can also assume there are no terms $z^{l_j^\prime}w^{k^\prime\alpha}w_j$ with $k^\prime>k$ and $l_j^\prime<l_j$ in the expression of $w_{j,1}$. This can be done because each blow-up centered at  $[1:0:\cdots:0]$ transforms terms like $z^{m}w^{k'\alpha}w_j$ in $w_{j,1}$ into terms like $z^{m+k'|\al|}w^{k'\alpha}w_j$. Thus, if $m\in \{0,\cdots, l_j-1\}$ is the smallest such $l_j'$, since $k'>k$, there exists $q\in \N$ such that $m+qk'|\al|>l_j+qk|\al|$, and thus, performing at most $q$ blow-ups centered at $[1:0:\cdots:0]$ we are done.

In any case, after holomorphic changes of coordinates and blow-ups, we can assume that $F$ has the form
\begin{equation}\label{total}
\left\{
\begin{aligned}
z_1 & = z - \frac{1}{\nu-1}z^{\nu} + o(z^\nu),\\
w_{j,1} & = \lambda_j w_j - a_jz^{l_j}w^{k\alpha}w_j - b_jz^{r_j}w_j + w_jo(z^{l_j}w^{k\alpha}, z^{r_j})+o(z^\mu),\ \ \ 1\le j\le n,
\end{aligned}
\right.
\end{equation}
where $l_j\geq 0$, $\mu>\nu$, and we use freely the Landau little/big-oh notation. For instance, the term $w_jo(z^{l_j}w^{k\alpha}, z^{r_j})$ denotes a holomorphic function of the form $w_jz^{l_j+1}w^{(k+1)\al}f_1(z,w)+w_jz^{r_j+1}f_2(z,w)$, for some holomorphic functions $f_1,f_2$.

As explained in the introduction, the condition $r_j=\nu-1$ for some $j$ was already studied in \cite{R:Quasi3} and the one-resonant part does not play any role. Therefore we assume that $F$ is \textit{degenerately dynamically separating}, namely, $r_j>\nu-1$. We also assume that
\[
l_1=\cdots=l_n=:l\leq \nu-1.
\]
The number $l\in \{0,\cdots, \nu-1\}$ is clearly an invariant under holomorphic changes of coordinates preserving the form \eqref{total}. Also, the condition $l_1=\cdots=l_n$ is invariant under blow-ups. However, the condition that $l\leq \nu-1$  is {\sl not} an invariant under blow-ups. Thus, starting from a given quasi-parabolic one-resonant germ, one has first to change coordinates to make it resonant up to high order, and then perform the least number  of blow-ups in order to get \eqref{total} and check then the condition on the $l_j$'s. As it will be clear from the proof of our theorem, if $l>\nu-1$ (or if some $l_j\neq l_k$) then the $z$-component enters too strongly in the behavior of the $w_{j,1}$ and we found no way to suitably control it.

In the given hypotheses, the number $A:=\sum_{j=1}^n \alpha_ja_j\lambda_j^{-1}$ is easily seen to be an invariant up to a scalar multiple under holomorphic changes of coordinates which preserve the form \eqref{total} (cf. \cite[Remark 3.2]{BZ:One}). We say that $F$ is \textit{non-degenerate with respect to $\{\lambda_1,\cdots,\lambda_n\}$} if $A\neq 0$.

Finally, in case $A\neq 0$, we say that $F$ is \textit{attracting with respect to $\{\lambda_1,\cdots,\lambda_n\}$} if $\Re(a_j\lambda_j^{-1}A^{-1})>0$ for all $1\le j\le n$. Once again, the condition of being attracting with respect to $\{\lambda_1,\cdots,\lambda_n\}$ is invariant under holomorphic changes of coordinates which preserve \eqref{total} (cf. \cite[Remark 5.2]{BZ:One}).

\section{Proof of Theorem \ref{T:Main}}

Let $F$ be as in Theorem \ref{T:Main}. Then, after holomorphic changes of coordinates and finitely many blow-ups if necessary, we can assume that $F$ is of the form
\begin{equation}\label{E:Good}
\left\{
\begin{aligned}
z_1 & = z - \frac{1}{\nu-1}z^{\nu} + o(z^\nu),\\
w_{j,1} & = \lambda_j w_j - a_jz^lw^{k\alpha}w_j - b_jz^\nu w_j +w_jo(z^{l}w^{k\alpha})+ w_jo(z^{\nu})+o(z^{\mu}),\ \ \ 1\le j\le n,
\end{aligned}
\right.
\end{equation}
where $\nu-1\geq l\geq 0$ and $\mu>\nu$ is arbitrarily large.

We first assume $l<\nu-1$. Let $0<\epsilon,  \epsilon',\beta<1$ (to be suitably chosen later), and let $\delta<1/k$. Let $\gamma$ be such that
\begin{equation}\label{defgamma}
\frac{k}{\nu-l-1}>\gamma>\frac{k}{\nu-l}.
\end{equation}
Moreover, if $l\geq 1$ let
\begin{equation}\label{defdelta}
0<\delta'<\frac{\delta}{2l(\nu-1)},
\end{equation}
otherwise let $0<\delta'<\frac{\delta}{2(\nu-1)}$.

For $a,b>0$ we set
$$V_{a,b}:=\{t\in \C:\ 0<|t|<a,\ |\arg(t)|<b\}.$$

Set $u:=w^\alpha=w_1^{\alpha_1}\cdots w_n^{\alpha_n}$. Let $\eta_1=1,\eta_2, \ldots, \eta_k$ be the roots of the equation $x^k=1$. For $t=1,\ldots, k$, define
\begin{equation}\label{E:B}
B_t:=\{(z,w)\in \C^{n+1}:\ |z|<|u|^\gamma,\ |w_j|<|u|^\beta,\ u\in \eta_tV_{\epsilon,\delta},\ z\in  V_{\epsilon^\prime,\delta^\prime}\}.
\end{equation}
The sets $B_t$'s are clearly disjoint because the projection $(z,w)\mapsto u$ maps them into disjoint sets.
We want to show that the $B_t$'s are open sets with $0\in\partial B_t$, that $F(B_t)\subset B_t$ and $F^n(p)\rightarrow O$ as $n\rightarrow \infty$ for $p\in B_t$, $t=1,\ldots, k$. This will prove Theorem \ref{T:Main}.

We will focus only on the set $B:=B_1$, the others being similar. First of all, since $(z,w)=(r^{2|\al|\gamma},r,\cdots,r)\in B$  for  $\R^+\ni r\to 0$,  the set $B$ is a non-empty open set and $O\in \partial B$.

Let $u_1=w_1^\alpha=w_{1,1}^{\alpha_1}\cdots w_{n,1}^{\alpha_n}$. Then we have
\begin{equation}\label{E:u1}
u_1=u(1-Az^lu^k-c z^\nu+ o(z^lu^k,z^\nu,\frac{z^\mu}{u})),
\end{equation}
with
\begin{equation}\label{E:A}
A:=\sum_{j=1}^n \alpha_ja_j\lambda_j^{-1}, \quad c:=\sum_{j=1}^n \alpha_jb_j\lambda_j^{-1}.
\end{equation}
Since $A$ is well-defined up to a non-zero scalar multiple, and $F$ is non-degenerate--thus $A\neq 0$--by re-scaling in the $w$ components if necessary, we can assume that $A=1/k$. Since $F$ is attracting with respect to $\{\lambda_1,\cdots,\lambda_n\}$, we have
\begin{equation}\label{E:Re}
\Re(a_j\lambda_j^{-1})>0.
\end{equation}
Therefore, up to choosing $\delta>0$ smaller if it is the case, we can also assume that
\begin{equation}\label{delta2}
\max_{j=1,\cdots, n} |\arg (a_j\lambda_j^{-1})|+2k\delta<\frac{\pi}{2}.
\end{equation}

Write $v=u^k$ and $v_1=u_1^k$. Now, given $(z,w)\in B$, by the very definition of $\gamma$ in \eqref{defgamma}, it is easy to get
\begin{equation}\label{little}
z^\nu=o(z^lv),\ \ \ z^lv=o(z^{\nu-1}).
\end{equation}
Moreover, note that when $(z,w)\in B$ then
\[
\frac{|z|^\mu}{|u|}=|z|^\nu \frac{|z|^{\mu-\nu}}{|u|}<|z|^\nu|u|^{(\mu-\nu)\gamma-1},
\]
and thus for $\mu$ large, we have that the term $o(\frac{z^\mu}{u})$ in the expression \eqref{E:u1} is in fact $o(z^lv)$ by \eqref{little}.

Therefore, from \eqref{E:u1} we obtain
\begin{equation}\label{E:u}
u_1=u(1-\frac{1}{k}z^lv+o(z^lv)),
\end{equation}
and, from this,
\begin{equation}\label{E:v}
v_1=v(1-z^lv+o(z^lv)).
\end{equation}
Therefore,
\begin{equation}\label{E:v1}
\frac{1}{v_1}=\frac{1}{v}+z^l+o(z^l).
\end{equation}

From (\ref{E:Good})  we have also
\begin{equation}\label{E:z}
\frac{1}{z_1^{\nu-1}}=\frac{1}{z^{\nu-1}}+1+O(z).
\end{equation}

For $R,r>0$ let us define
$$U_{R,r}:=\{t\in \C:\ |t|>R,\ |\arg(t)|<r\}.$$
Note that $x\in V_{a,r}$ if and only if $1/x^k\in U_{a^{-k},kr}$.

Now we want to show that if $(z,w)\in B$, then  $z_1\in V_{\epsilon^\prime,\delta^\prime}$ and  $u_1\in  V_{\epsilon,\delta}$. Set $R=\epsilon^{-k}$ and $R'=\epsilon'^{-(\nu-1)}$. By what we said above, this is equivalent to showing that $1/z_1^{\nu-1}\in U_{R^\prime,(\nu-1)\delta^\prime}$ and  $1/v_1\in U_{R,k\delta}$.

Note that, if $\tau\in \C$ is such that $|\tau|<\tan((\nu-1)\delta')$, then $U_{R^\prime,(\nu-1)\delta^\prime}+1+\tau\subset U_{R^\prime,(\nu-1)\delta^\prime}$. From (\ref{E:z}), let $C>0$ be such that $|\frac{1}{z_1^{\nu-1}}-\frac{1}{z^{\nu-1}}-1|\leq C$ for $(z,w)\in B$. Hence if $\epsilon'>0$ is such that $\epsilon'C<\tan( (\nu-1)\delta')$, it follows that
\begin{equation}\label{z-inside}
1/z^{\nu-1}\in U_{R^\prime,\delta^\prime}\Rightarrow 1/z_1^{\nu-1}\in U_{R^\prime,(\nu-1)\delta^\prime}.
\end{equation}

As for $1/v_1$, if $l=0$ it is clear. If $l\geq 1$, note that since $z\in V_{\epsilon^\prime,\delta^\prime}$, it follows that $z^l\in V_{\epsilon^\prime,l\delta^\prime}\subset V_{\epsilon^\prime,\delta}$ by \eqref{defdelta}. Note also that if $x\in U_{R,k\delta}$ then $x+V_{\epsilon^\prime,\delta}\subset U_{R,k\delta}$.  Hence, from \eqref{E:v1}, given $(z,w)\in B$, choosing $\epsilon'$ smaller if necessary,
it follows that  $1/v_1\in U_{R,k\delta}$, as claimed.

Now, we want to show that $|w_{j,1}|<|u_1|^\beta$ for $j=1,\cdots, n$ for $(z,w)\in B$. From (\ref{E:Good}) and (\ref{E:u}), taking into account \eqref{little}, we  have
$$\begin{aligned}
\frac{w_{j,1}}{u_1^{\beta}} & =\frac{w_j}{u^{\beta}}\frac{\lambda_j-a_jz^lv+o(z^lv)}{1-\frac{\beta}{k} z^lv+o(z^lv)}\\
& =\frac{w_j\lambda_j}{u^{\beta}}\left(1-(\frac{a_j}{\lambda_j}-\frac{\beta}{k})z^lv+o(z^lv)\right).
\end{aligned}$$
Choosing $\beta<k \min_{j=1,\cdots, n}\{\Re a_j\lambda_j^{-1}\}$, by  \eqref{E:Re}, we obtain that $\Re(\frac{a_j}{\lambda_j}-\frac{\beta}{k})>0$ for $j=1,\cdots, n$. Also, $z^l\in V_{\epsilon^\prime,\delta}$ while $v=u^k\in V_{\epsilon, k\delta}$ and by \eqref{delta2} it follows that
\[
\Re \left((\frac{a_j}{\lambda_j}-\frac{\beta}{k})z^lv\right)>0.
\]
From this, for $\epsilon, \epsilon'<<1$, we have $|w_{j,1}|<|u_1|^\beta$ for $j=1,\cdots, n$ as claimed.

Finally, we want to show that $|z_1|<|u_1|^\gamma$ if $(z,w)\in B$.

From (\ref{E:Good}) and \eqref{E:u} (and \eqref{little}), we have
\begin{equation}\label{z-1}\begin{aligned}
\frac{z_1}{u_1^{\gamma}} & =\frac{z}{u^{\gamma}}\frac{1-\frac{1}{\nu-1}z^{\nu-1}+o(z^{l}v)}{1-\frac{\gamma}{k} z^lv+o(z^lv)}\\
& =\frac{z}{u^{\gamma}}(1-\frac{1}{\nu-1}z^{\nu-1}+\frac{\gamma}{k}z^lv+ o(z^lv)).
\end{aligned}
\end{equation}

Now, since $(z,w)\in B$ and hence $|z|<|u|^\gamma$, we can write $|z|=|u|^{\theta}$ for some $\theta>\gamma$. If $\theta<\frac{k}{\nu-l-1}$ then $z^lv=o(z^{\nu-1})$ as in \eqref{little}. In this case, from \eqref{z-1} we have that
\[
\frac{|z_1|}{|u_1|^{\gamma}} =\frac{|z|}{|u|^{\gamma}}|1-\frac{1}{\nu-1}z^{\nu-1}+o(z^{\nu-1})|.
\]
Since $\Re z^{\nu-1}>0$, being $z^{\nu-1}\in V_{\epsilon^\prime,(\nu-1)\delta^\prime}\subset V_{\epsilon^\prime,\delta}$ by \eqref{defdelta}, if $\epsilon,\epsilon'<<1$ we have that $|1-\frac{1}{\nu-1}z^{\nu-1}+o(z^{\nu-1})|<1$ and hence
\[
|1-\frac{1}{\nu-1}z^{\nu-1}+o(z^{\nu-1})|<\frac{|z|}{|u|^{\gamma}}<1,
\]
as needed.

On the other hand, if $|z|=|u|^{\theta}$ for some $\theta\geq \frac{k}{\nu-l-1}>\gamma$, then from \eqref{z-1} we have
\[
\frac{|z_1|}{|u_1|^{\gamma}}=|u|^{\theta-\gamma}(1+o(|u|^{\theta-\gamma})),
\]
hence if $\epsilon,\epsilon'<<1$ it follows that $|u|^{\theta-\gamma}(1+o(|u|^{\theta-\gamma}))<1$ and we are done.

We proved therefore that $F(B)\subset B$. Hence, applying inductively (\ref{E:z}) and (\ref{E:v1}) to a point $(z,w)\in B$, we get
\begin{equation}\label{E:zn}
\frac{1}{|z_n|^{\nu-1}}\sim n\ \Rightarrow\ |z_n|\sim (\frac{1}{n})^{\frac{1}{\nu-1}},
\end{equation}
and
\begin{equation}\label{E:un}
\frac{1}{|v_n|}\sim \sum_{j=1}^n |z_j|^l\sim \sum_{j=1}^n (\frac{1}{j})^{\frac{l}{\nu-1}}\sim n^{1-\frac{l}{\nu-1}}\ \Rightarrow\ |u_n|\sim (\frac{1}{n})^{\frac{1-\frac{l}{\nu-1}}{k}}.
\end{equation}

Now we consider the case $l=\nu-1$. Note that in this case $l\geq 1$. We retain the previously introduced notations. The proof goes similarly to the previous case, except that we have to choose
\begin{equation}\label{defg2}
0<\delta'<\frac{\delta}{2l^2},\quad \gamma<\frac{k}{\sqrt{1+\tan 1}}
\end{equation}
and define for $t=1,\ldots, k$,
\begin{equation}\label{E:B1}
B_t=\{(z,w)\in \C^{n+1}:\ |u|^\gamma\log |z|<-\frac{1}{l},\ |w_j|<|u|^\beta,\ u\in \eta_t V_{\epsilon,\delta},\ z\in V_{\epsilon^\prime,\delta^\prime}\}.
\end{equation}
As before, we concentrate on the case $B:=B_1$. Now, given $(z,w)\in B$ we have
$$\frac{1}{|z|^\nu}=\frac{1}{|z|^l}\frac{1}{|z|}>\frac{1}{|z|^l}(-l\log |z|)^{k/\gamma}>\frac{1}{|z|^l}|u|^{-k}=\frac{1}{|z|^l|v|},$$
and thus
\begin{equation}\label{little2}
z^\nu=o(z^lv).
\end{equation}
Using \eqref{little2} instead of \eqref{little}, starting from $(z,w)\in B$ and arguing as before we obtain immediately that $z_1\in V_{\epsilon^\prime,\delta^\prime}$, $u_1\in V_{\epsilon,\delta}$ and $|w_{j,1}|<|u_1|^\beta$.

We are only left to show that $|u_1|^\gamma\log |z_1|<-1/l$. First of all, note that if $z\in V_{\epsilon^\prime,\delta^\prime}$ hence $z^l\in V_{\epsilon^\prime,l\delta^\prime}\subset V_{\epsilon,\delta}$, hence
\begin{equation}\label{estz}
|z|\leq \sqrt{1+\tan(\delta')}\Re z, \quad |z|^l\leq \sqrt{1+\tan(\delta)}\Re (z^l).
\end{equation}
In particular, $|z|\sim \Re z$ and using \eqref{E:u} and \eqref{E:Good} (with $l=\nu-1$), we can write
\begin{equation}\label{E:log}
\begin{aligned}
|u_1|^\gamma\log |z_1| & =|u|^\gamma|1-\frac{1}{k}z^lv+o(z^lv)|^\gamma(\log |z|+\log |1-\frac{1}{l}z^l+o(z^l)|)\\
&=|u|^\gamma|1-\frac{1}{k}z^lv+o(z^lv)|^\gamma(\log |z|+\Re \log (1-\frac{1}{l}z^l+o(z^l)))\\
& =|u|^\gamma(\log |z|)|1-\frac{\gamma}{k} z^lv|(1-\frac{\Re z^l}{l\log |z|})+o(z^lv,\frac{z^l}{\log |z|})\\&=:|u|^\gamma(\log z) R(z,w).
\end{aligned}
\end{equation}

Now, since $(z,w)\in B$ and hence $\ds -\frac{1}{l\log |z|}<|u|^\gamma$, we have $\ds -\frac{1}{l\log |z|}=|u|^{\theta}$ for some $\theta>\gamma$.

If $\theta<k$, then $R(z,w)=1+\Re (z^l)|v|^{\theta/k}+o(z^lv^{\theta/k})>1$ because $\Re z^l>0$ since $z^l\in V_{\epsilon', l\delta'}\subset V_{\epsilon, \delta}$, and thus by \eqref{E:log} we get $|u_1|^\gamma\log |z_1|<|u|^\gamma\log |z|<-1/l$.

If $\theta=k$, then $R(z,w)=|1-\frac{\gamma}{k}z^lv|(1+|v|\Re z^l)+o(z^lv)$. We want to show that $R(z,w)\geq 1$. To this aim, set $\eta:=(1+\tan 1)^{-1/2}-\gamma/k$ and note that $\eta>0$ by \eqref{defg2}.  By \eqref{estz} and since the function $(0,1)\ni \delta\mapsto (1+\tan \delta)^{-1/2}$ is decreasing, we have
\begin{equation*}
\begin{split}
|1-\frac{\gamma}{k}z^lv|^2(1+|v|\Re z^l)^2&=(1-2\frac{\gamma}{k}\Re (z^lv))(1+2|v|\Re z^l)+o(vz^l)\\
&\geq (1-2\frac{\gamma}{k}|z^lv|)(1+\frac{2|z^lv|}{\sqrt{1+\tan \delta}})+o(vz^l)\\
&=1+2\left(\frac{1}{\sqrt{1+\tan \delta}}-\frac{\gamma}{k}\right)|z^lv|+o(vz^l)\\
&\geq 1+2\eta|z^lv|+o(vz^l),
\end{split}
\end{equation*}
from which it follows immediately that $R(z,w)\geq 1$ and again we get $|u_1|^\gamma\log |z_1|<-1/l$.

If $\theta>k$, then $|u_1|^\gamma\log |z_1|=|u|^\gamma\log |z| R(z,w)=-1/l|u|^{\gamma-\theta}R(z,w)<-1/l$, since $R(z,w)$ is close to 1, $|u|$ is small and $\gamma-\theta<\gamma-k<0$.

We proved therefore that $F(B)\subset B$. Hence, applying inductively (\ref{E:z}) and (\ref{E:v1}) to a point $(z,w)\in B$, we get \eqref{E:zn} and
\begin{equation}\label{E:un1}
\frac{1}{|v_n|}\sim \sum_{j=1}^n |z_j|^{\nu-1}\sim \sum_{j=1}^n \frac{1}{j}\sim \log n\ \Rightarrow\ |u_n|\sim (\frac{1}{\log n})^{\frac{1}{k}}.
\end{equation}

From \eqref{E:zn}, \eqref{E:un} and \eqref{E:un1}, it follows that for any $(z,w)\in B$ we have $F^n(z,w)\rightarrow O$ as $n\rightarrow \infty$. This completes the proof of Theorem \ref{T:Main}.

\begin{remark}\label{R:l}
Let $\varrho_1=1,\varrho_2, \ldots, \varrho_{l}$ be the roots of the equation $x^l=1$. For $t=1,\ldots, k$ and $s=1,\ldots, \nu-1$ define
\[
B_{t,s}=\{(z,w)\in \C^{n+1}:\ |u|^\gamma\log |z|<-\frac{1}{l},\ |w_j|<|u|^\beta,\ u\in \eta_t V_{\epsilon,\delta},\ z\in \varrho_s V_{\epsilon^\prime,\delta^\prime}\}.
\]
Those $B_{t,s}$'s are disjoint open sets with $0\in \partial B_{t,s}$. In case $l=\nu-1$, since $z^l=z^{\nu-1}$ belongs to $V_{\epsilon^\prime,l\delta^\prime}$ whenever $z\in \varrho_s V_{\epsilon^\prime,\delta^\prime}$,  the previous proof shows  that $F(B_{t,s})\subset B_{t,s}$ and $F^{n}(p)\to 0$ as $n\to \infty$ for $p\in B_{t,s}$. Thus, in case $l=\nu-1$, $F$ has (at least) $(\nu-1)k$ disjoint basins of attraction. Similar arguments can be carried out for the case $l=0$.

In case $0<l<\nu-1$ the previous argument fails because we cannot control $z^l$ since it does not stay in one ``petal" and hence, from \eqref{E:v1}, we cannot infer anything about the behavior of $v_1$.
\end{remark}

\end{document}